\documentclass[11pt, letterpaper]{article}

\usepackage{amsmath,amsthm}
\usepackage{enumerate}
\usepackage{multirow}
\usepackage{geometry}
\usepackage[T1]{fontenc}
\usepackage{array}
\usepackage{makecell}
\newtheorem{theorem}{Theorem}

\title{Relative error due to a single bit-flip in floating-point arithmetic}
\author{Bradley R. Lowery}

\begin{document}
\maketitle

We consider the error due to a single bit-flip in a floating point number.  We
assume IEEE 754 double precision arithmetic, which encodes binary floating
point numbers in a 64-bit word~\cite{ieee}. We assume that the bit-flip
happens randomly so it has equi-probability (1/64) to hit any of the 64 bits.
Since we want to mitigate the assumption on our 
initial floating-point number, we assume that it is uniformly picked among all
normalized number. With this framework, we can summarize our findings as follows.
\begin{center}\begin{minipage}{13cm}
The probability for a single bit flip to cause a relative error less than $10^{-11}$ in a normalized floating-point number is above 25\%;
The probability for a single bit flip to cause a relative error less than $10^{-6}$ in a normalized floating-point number is above 50\%; Etc.
\end{minipage}\end{center}

Similar studies could be done for binary floating point numbers in a 32-bit
word, or complex numbers; (so, s,c,d,z in LAPACK or BLAS notation); multiple
bit-flips; or burst of bit-flips (that is multiple bit-flips which happen
contiguously in the binary word).

A 64-bit word consists of 3 fields: a sign bit ($s$), a biased 11 bit exponent
($e$), and a 52 bit fraction ($f$).  Excluding exceptions (to be explained) the
encoded number is \[ (-1)^s(1.f)_2 2^{e-1023}. \] There are three exceptions:
not a number (NAN), infinity (INF), and denormalized numbers.  NANs are
represented by the biased exponent of all ones and any nonzero fraction (sign
does not matter).  INFs are represented by the biased exponent of all ones and
the zero fraction (positive or negative infinity is possible).  Denormalized
numbers have a biased exponent of all zeros.  In this case the number
represented is $(-1)^s(0.f)_2 2^{-1022}$.  Zero is a special case of a
denormalized number where the fraction is zero. We refer the reader to
Overton~\cite{overton} for more details on IEEE arithmetic.

We only make one mild assumption on the initial floating-point number. For
example, we assume it is a normalized number and then derive probabilities from
that assumption. In that case, we assume that our initial floating-point number
is randomly picked among all normalized number.  In some cases, one could
derive more information if one assumes more on the initial floating-point
number. For example, if the number is $2.0e+00$, then the probability that a
bit-flip changes it to INF is zero. That would be a different study.

Let $x$ be the original 64-bit word and $x'$ be the 64-bit word after a single
bit-flip between $x$ and $x'$ occurred.  We will also denote the fields of $x'$
as $s'$, $e'$, and $f'$ when needed.  Also, $e(k)$ is the $k$-th entry of the
exponent field, and similar notation for other fields.

We begin by calculating the probability that 
$x'$ is normalized, NAN, INF, or denormalized, 
given that $x$ is normalized, NAN, INF, or denormalized.
Table~\ref{t.changetype} summarizes the results.
The narrative below explains where the results come from.

\begin{enumerate}
	\item Let $x$ be a normalized number. Therefore, the exponent is anything other than
	all ones or all zeros. 
	For $x' = \textmd{NAN}$ or $x' = \textmd{INF}$ the
	exponent must be all ones. There are 11 exponents corresponding to normalized numbers
	that are susceptible to becoming all ones due to a bit-flip 
	(bit strings of all ones except 1 zero).  The probability 
	that the zero bit in the exponent is flipped is $1/64$.  The result will be an INF
	only if $f = 0$.  Therefore we have
	\[ \Pr( \textmd{normalized to NAN } ) 
		= \left(\frac{11}{2^{11}-2}\right)\left(\frac{1}{64}\right)\left( 1 - \frac{1}{2^{52}} \right) 
		\approx 8.4005\times 10^{-5}
	\]
	and
	\[ \Pr( \textmd{normalized to INF } ) 
		= \left(\frac{11}{2^{11}-2}\right)\left(\frac{1}{64}\right)\left( \frac{1}{2^{52}} \right)
		\approx 1.8653\times 10^{-20}
	. \]
	For $x'$ to be a denormalized number the exponent must be all zeros. There are 
	11 exponents corresponding to normalized numbers
	that are susceptible to becoming all zeros due to a bit-flip 
	(bit strings of all zeros except 1 one).  The probability 
	that the one bit in the exponent is flipped is $1/64$.  For all factions the result
	will be a denormalized number.  Therefore we have 
	\[ \Pr( \textmd{normalized to denormalized } ) 
		= \left(\frac{11}{2^{11}-2}\right)\left(\frac{1}{64}\right)
		\approx 8.4005\times 10^{-5}
	. \]
	Finally, we consider the probability that $x'$ is a normalized number.
	If the exponent of $x$ is not one of the 22 susceptible values then $x'$ will be a
	normalized number. $x'$ will also be a normalized number if the exponent is one of the
	susceptible exponents but the bit-flip is not in the corresponding entry. 
	\[ \Pr( \textmd{normalized to normalized } ) 
		= 1 - \left(\frac{22}{2^{11}-2}\right)\left(\frac{1}{64}\right)
		\approx 9.9983\times 10 ^{-1}
	. \]
	\item Let $x = \textmd{NAN}$. Therefore, $e = (1\dots1)_2$ and $f \ne 0$.
	$x'$ will only be a normalized number if a bit in the exponent is flipped. Therefore,
	\[ \Pr( \textmd{NAN to normalized} ) = \frac{11}{64}
		\approx 1.7188\times 10^{-1}
	. \]
	There are 52 fractions that are susceptible to becoming all zeros (bit strings
	of all zeros except 1 one).  The probability that the 1 bit is flipped is $1/64$. 
	If this occurs the result will be an INF, therefore 
	\[ \Pr( \textmd{NAN to INF} ) = \left( \frac{52}{2^{52}} \right)
		\left( \frac{1}{64} \right)
		\approx 1.8041\times 10^{-16}
	. \]
	A single bit-flip can not produce a denormalized number, therefore
	\[ \Pr( \textmd{NAN to denormalized} ) = 0. \]
	The remaining cases will result in a NAN, therefore
	\[ \Pr( \textmd{NAN to NAN} ) = 1 - \frac{11}{64}
		- \left( \frac{52}{2^{52}} \right)
		\left( \frac{1}{64} \right)
		\approx 8.2812\times 10^{-1}
	. \]

	\item Let $x = \textmd{INF}$.  Therefore, $e = (1\dots1)_2$ and $f = (0\dots0)_2$.
	If the bit-flip occurs in $f$, then $x' = \textmd{NAN}$.
	If the bit-flip occurs in $e$, then $x'$ is a normalized number.
	If the bit-flip occurs in $s$, then $x' = \pm\textmd{INF}$.
	No bit-flip will cause a denormalized number.  Therefore, we have the following probabilities:
	\begin{gather*}
		\Pr( \textmd{INF to NAN} ) = 52/64 		\approx 8.1250\times 10^{-1} , \\
		\Pr( \textmd{INF to normalized} ) = 11/64	\approx 1.7188\times 10^{-1} , \\
		\Pr( \textmd{INF to INF} ) = 1/64		\approx 1.5625\times 10^{-2} , \\
		\Pr( \textmd{INF to denormalized} ) = 0. \\
	\end{gather*}

	\item Let $x$ be a denormalized number. Therefore, $e = (0\dots0)_2$. If the bit-flip occurs
	in either $s$ or $f$ then $x'$ will also be a denormalized number.  If the bit-flip occurs in
	$e$, then $x'$ will be a normalized number.  No bit-flip will cause an NAN or INF.
	Therefore, we have the following probabilities:
	\begin{gather*}
		\Pr( \textmd{denormalized to denormalized} ) = 53/64	\approx 8.2812\times 10^{-1} , \\
		\Pr( \textmd{denormalized to normalized} ) = 11/64		\approx 1.7188\times 10^{-1}, \\
		\Pr( \textmd{denormalized to NAN} ) = 0, \\
		\Pr( \textmd{denormalized to INF} ) = 0. \\
	\end{gather*}
\end{enumerate}

\begin{table}
	\centering
	\begin{tabular}{ |c | c | c | c | c | }
		\hline
		\diaghead(5,-1){Diag ColumnHead II}{$x$}{$x'$}	&	Normalized	&	Denormalized	&	NAN	&	INF \\ \hline
		Normalized	
		& $9.9983\times 10 ^{-1}$
		& $8.4005\times 10^{-5}$
		& $8.4005\times 10^{-5}$
		& $1.8653\times 10^{-20}$
		\\ \hline
		Denormalized 
		& $1.7188\times 10^{-1}$	
		& $8.2812\times 10^{-1}$
		& 0
		& 0
		\\ \hline
		NAN
		& $1.7188\times 10^{-1}$
		& 0
		& $8.2812\times 10^{-1}$
		& $1.8041\times 10^{-16}$
		\\ \hline
		INF
		& $1.7188\times 10^{-1}$ 
		& 0
		& $8.1250\times 10^{-1}$
		& $1.5625\times 10^{-2}$
		\\ \hline
	\end{tabular}
	\caption{Probability of a 64 bit word changing type due to a single bit-flip.}
	\label{t.changetype}
\end{table}

\section*{Relative error}

\subsection*{Normalized Numbers}
When $x$ and $x'$ are normalized or denormalized numbers we would like to calculate
the relative error in $x'$ as well as the probability that the relative error is below a given threshold.

Let $x$ and $x'$ be normalized numbers.
Then we have the following cases.
% for the bit-flip occurring in each field.
\begin{enumerate}

\item If the flipped bit is in $s$, then $x' = -x$ and
\[ \frac{| x - x' |}{| x |} = 2. \]
\item If the flipped bit is the $k$-th entry of $f$, then
\[ | x - x' | = (0.0\dots010\dots0)_2 2^{e-1023} = (1.0)_2 (2^{-k}) 2^{e-1023}. \]
Hence,
\[ \frac{| x - x' |}{| x |} = \frac{(1.0)_2}{(1.f)_2} (2^{-k}). \]
Since $ 1/2 < 1 / (1.f)_2 \le 1 $, we have  
\[ 2^{-k-1} < \frac{| x - x' |}{| x |} \le 2^{-k}. \]
\item If the flipped bit is in $k$-th entry of $e$, then
\begin{align*}
	| x - x' | 	&= | (1.f)_2 (2^{-1023})(2^e - 2^{e'}) | \\
				&= | (1.f)_2(2^{-1023})(2^e) | | 1 - 2^{e' - e} | \\
				&= | x | | 1 - 2^{e' - e} |. 
\end{align*}
	\begin{enumerate}
		\item If $e(k) = 0$, hence $e' > e$, then $e' - e = 2^{11-k}$.
		\item If $e(k) = 1$, hence $e' < e$, then $e' - e = -2^{11-k}$.
	\end{enumerate}
Therefore,
\begin{equation*}
\frac{| x - x' |}{| x |} =
	\begin{cases}
		2^{2^{11-k}} - 1 & \textmd{if} \; x(k) = 0, \\
		1 - 2^{-2^{11-k}} & \textmd{if} \; x(k) = 1. 
	\end{cases}
\end{equation*}
Note that $2^{2^{11-k}} - 1 \ge 1$  
	and $ 1/2 \le 1 - 2^{-2^{11-k}} < 1$ for $k = 1,\dots,11$.
\end{enumerate}

Let $x$ be a normalized number and $x'$ be a denormalized number.  This 
can only happen if the bit-flip occurs in the exponent. The error is
\begin{align*}
	| (-1)^s(1.f)_2 2^{e-1023} - (-1)^s(0.f)_2 2^{-1022} |
		&= | x | | 1 - \frac{(0.f)_2}{(1.f)_2} 2^{1-e} |
\end{align*}
Since, $0 \le (0.f)_2 / (1.f)_2 < 1/2$ and $e = 2^{11-k}$  
\[ \frac{1}{2} \le 1 - 2^{-e} < \frac{| x - x' |}{| x |} \le 1. \]
The following theorem summarizes the results.

\begin{theorem}\label{thm.normalized}
Let $x$ is a normalized number.  The relative error due to a single bit-flip is
\begin{alignat}{2}
		| x - x' | / | x | = 2 
				&\quad \textmd{if bit-flip in $s$,}  \label{eq.n1}\\
		2^{-k-1} < | x - x' | / | x | \le 2^{-k} 
				&\quad\textmd{if bit-flip in $f$,} \label{eq.n2} \\
		| x - x' | / | x | = 2^{2^{11-k}} - 1 \ge 1 
				&\quad\textmd{if bit-flip in $e$ and $e(k) = 0$,} \label{eq.n3}\\
		1/2 \le | x - x' | / | x | = 1 - 2^{-2^{11-k}} < 1
				&\quad\textmd{if bit-flip in $e$, $e(k) = 1$, and $e' \ne 0$} \label{eq.n4}\\
		1/2  \le 1 - 2^{-2^{11-k}} < | x - x' | / | x |  \le 1
				&\quad\textmd{if bit-flip in $e$ and $e' = 0$,} \label{eq.n5}
\end{alignat}
where $k$ is the location of the bit-flip in the given field.
\end{theorem}

Using these results,
assuming $x$ is a random normalized number,
we can calculate the probability for the relative error to be in the following intervals
\begin{enumerate}
	\item $\Pr(| x - x' | / | x | \ge 1 ) = 13/128 + (11 / (2046)2^{-52}(1/64) \approx 0.10156$. 
	This is mainly due to the probability of a bit-flip in $s$ and in $e$ (flipping 0 to 1).
	See \eqref{eq.n1} and \eqref{eq.n3}.  The probability of each of these cases is $1/64$ and 
	$11/128$, respectively.
	There is also a possibility for the relative error
	to be 1 if $e' = (0\dots0)_2$ and $f = 0$, hence $x' = 0$.  This is a special case of \eqref{eq.n5}.
	This will occur with a probability of $(11/2046)(1/64)2^{-52}$. 
	\item $\Pr( 1/2 < | x - x' | / | x | < 1 ) = 11/128 - ((2^{10}-2)/2046)(1/64) 
		- (11 / (2046))2^{-52}(1/64) \approx 0.078133$. 
	This is mainly due to the probability of a bit-flip in $e$ (flipping 1 to 0), which occurs
	with probability $11/128$. See \eqref{eq.n4} and \eqref{eq.n5}.  
	However, these include the special case when the error is $1/2$ and $1$, which we do not include.
	The error will be $1/2$ from this case if the first bit of the exponent is flipped, but the
	resulting number must still be normalized.  This occurs with probability $(2^{10}-2)/2046)(1/64)$.  
	The error will be 1 for this case if $x' = 0$, which occurs with probability 
	$(11/2046)(1/64)2^{-52}$.
	\item $\Pr( | x - x' | / | x | \le 1/2 ) = 52/64 + (2^{10}-2)(2046)(1/64) \approx 0.82030$. 
	The relative error will always be less
		than or equal to $1/2$ if the bit-flip occurs in $f$ (see \eqref{eq.n2}),
		which occurs with probability $52/64$.
		Including the special case from \eqref{eq.n4} such that the relative error is $1/2$ gives
		the stated probability.
	\item $\Pr( | x - x' | / | x | \le 2^{-i} ) = (53 - i)/64$ for $i = 2,\dots,52$. 
		The relative error will always be less than or equal to $2^{-i}$ if the bit-flip
		occurs in the $j$-th entry of $f$ ($j \ge i$). See \eqref{eq.n2}.
		 There are $53-i$ entries of $f$ greater than or equal to $i$. 
\end{enumerate}

We can bound the probability that the error is less than some tolerance using the fourth case above. 
Below are the bounds on the probability for some error tolerances.

\begin{align*}
	\Pr( | x - x' | / | x | \le 10^{-1} )  &< 0.78125 &
		\Pr( | x - x' | / | x | \le 10^{-9} )  &< 0.37500 \\ 
	\Pr( | x - x' | / | x | \le 10^{-2} )  &< 0.73438 & 
		\Pr( | x - x' | / | x | \le 10^{-10} ) &< 0.31250 \\ 
	\Pr( | x - x' | / | x | \le 10^{-3} )  &< 0.68750 & 
		\Pr( | x - x' | / | x | \le 10^{-11} ) &< 0.26562 \\ 
	\Pr( | x - x' | / | x | \le 10^{-4} )  &< 0.62500 & 
		\Pr( | x - x' | / | x | \le 10^{-12} ) &< 0.21875 \\ 
	\Pr( | x - x' | / | x | \le 10^{-5} )  &< 0.57812 & 
		\Pr( | x - x' | / | x | \le 10^{-13} ) &< 0.15625 \\ 
	\Pr( | x - x' | / | x | \le 10^{-6} )  &< 0.53125 & 
		\Pr( | x - x' | / | x | \le 10^{-14} ) &< 0.10938 \\ 
	\Pr( | x - x' | / | x | \le 10^{-7} )  &< 0.46875 & 
		\Pr( | x - x' | / | x | \le 10^{-15} ) &< 0.06250 \\
	\Pr( | x - x' | / | x | \le 10^{-8} )  &< 0.42188 & &
\end{align*} 

\subsection*{Denormalized Numbers}
Let $x$ be a denormalized number. 
If the bit-flip occurs in $s$, then $x'$ is also a denormalized number and the
error is the same as in the normalized case: 
\[ \frac{| x - x' |}{| x |} = 2. \]

If the bit-flip occurs in $k$-th entry of $f$, then $x'$ will again 
be a denormalized number and 
\[ | x - x' | = (0.0\dots010\dots0)_2 2^{-1022} = (1.0)_2 (2^{-k}) 2^{-1022}. \]
Hence,
\[ \frac{| x - x' |}{| x |} = \frac{(1.0)_2}{(0.f)_2} (2^{-k}). \]
Since $ 1 < 1 / (0.f)_2 \le 2^{52} $, we have  
\[ 2^{-k} < \frac{| x - x' |}{| x |} \le 2^{52-k}. \]
To improve these bounds we must be more specific on $f$.
Let the first nonzero entry in $f$ be the $t$-th entry.  
Then $ 2^{t-1} < 1 / (0.f)_2 \le 2^{t} $ and
\[ 2^{t-k-1} < \frac{| x - x' |}{| x |} \le 2^{t-k}. \]
Finally, if the bit-flip occurs in then $k$-th entry of  $e$, 
then $x'$ is a normalized number and 
\begin{align*}
	| x - x' | 	&= | (0.f)_2 2^{-1022} - (1.f)_2 2^{e'-1023} | \\
				&= | 2^{-1022}( (0.f)_2 - 2^{-1}(1.f)_2 2^{e'} | \\
				&= | 2^{-1022}( (0.f)_2 - (0.f)_2 2^{e'} | \\
				&= | 2^{-1022}(0.f)_2 ( 1 -  2^{e'} | \\
				&= | x | | 1 -  2^{e'} | \\
				&= | x | | 1 -  2^{2^{11-k}} |. \\
\end{align*}
For the last line we use $e' = 2^{11-k}$.  This is the same error as in the
normalized case when the bit-flip in the exponent is from 0 to 1. The following
theorem summarizes the results.  
\begin{theorem}
Let $x$ is a denormalized number.  The relative error due to a single bit-flip is
\begin{alignat}{2}
		| x - x' | / | x | = 2 
				&\quad \textmd{if bit-flip in $s$,} \\
		2^{t-k-1} < | x - x' | / | x |  \le 2^{t-k}
				&\quad \textmd{if bit-flip in $f$,} \\
		| x - x' | / | x |  = 2^{2^{11-k}} - 1 \ge 1 
				&\quad \textmd{if bit-flip in $e$,} 
\end{alignat}
where $k$ is the location of the bit-flip in the given field
and $t$ is the location of the first nonzero entry in $f$.

\end{theorem}

\bibliographystyle{plain}
\bibliography{biblio}

\begin{thebibliography}{1}

\bibitem{ieee}
{IEEE} standard for binary floating-point arithmetic.
\newblock {\em ANSI/IEEE Std 754-1985}, 1985.

\bibitem{overton}
M.L. Overton, Society for Industrial, and Applied Mathematics.
\newblock {\em Numerical Computing with IEEE Floating Point Arithmetic:
  Including One Theorem, One Rule of Thumb, and One Hundred and One Exercises}.
\newblock Society for Industrial and Applied Mathematics, 2001.

\end{thebibliography}

\end{document}